\documentclass{amsart}

\usepackage{amsmath}
\usepackage{amsthm}
\usepackage{amssymb}

\usepackage{comment}

\usepackage{hyperref}

\newtheorem{theorem}{Theorem}[section]
\newtheorem{definition}[theorem]{Definition}
\newtheorem{lemma}[theorem]{Lemma}
\newtheorem{corollary}[theorem]{Corollary}
\newtheorem{proposition}[theorem]{Proposition}
\newtheorem{remark}[theorem]{Remark}
\newtheorem{hypothesis}[theorem]{Hypothesis {\bf H.}\hspace*{-0.6ex}}

\newcommand{\R}{{\mathbb R}}
\newcommand{\N}{{\mathbb N}}
\newcommand{\Z}{{\mathbb Z}}
\newcommand{\C}{{\mathbb C}}

\newcommand{\nn}{\nonumber}
\newcommand{\be}{\begin{equation}}
\newcommand{\ee}{\end{equation}}

\newcommand{\ol}{\overline}
\newcommand{\ti}{\tilde}

\newcommand{\spr}[2]{\langle #1 , #2 \rangle}

\newcommand{\E}{\mathrm{e}}

\newcommand{\sgn}{\mathrm{sgn}}
\newcommand{\tr}{\mathrm{tr}}
\newcommand{\im}{\mathrm{Im}}

\DeclareMathOperator{\Ran}{Ran}

\newcommand{\abs}[1]{\lvert#1 \rvert}

\newcommand{\eps}{\varepsilon}

\newcommand{\sig}{\sigma}
\newcommand{\lam}{\lambda}
\newcommand{\gam}{\gamma}

\numberwithin{equation}{section}

\begin{document}

\title{On Perturbations of Quasiperiodic Schr\"odinger Operators}

\author{Helge Kr\"uger}
\address{Department of Mathematics\\ Rice University\\ Houston\\ TX 77005\\ USA}
\email{\href{mailto:helge.krueger@rice.edu}{helge.krueger@rice.edu}}
\urladdr{\href{http://www.mat.univie.ac.at/~helge/}{http://www.mat.univie.ac.at/\~{}helge/}}



\begin{abstract}
Using relative oscillation theory and the reducibility result of Eliasson,
we study perturbations of quasiperiodic Schr\"odinger operators. In particular,
we derive relative oscillation criteria and eigenvalue asymptotics for
critical potentials.
\end{abstract}

\maketitle

\section{Introduction}

We will be interested in generalizing classical perturbation result
of eigenvalues to quasiperiodic operators. We first overview
the classical results of interest. Most (if not all)
of our results will be parallel to these.
For this introduction let $H$ be a self-adjoint realization
of
\be
H = - \frac{d^2}{dx^2} + q(x)
\ee
on $L^2(1,\infty)$ with $q(x) \to 0$ as $x\to\infty$ and $q$ bounded.
A classical result of Weyl now tells us, that the essential spectrum
of $H$, is equal to the one of $-\frac{d^2}{dx^2}$, hence 
$\sig_{ess}(H) = [0, \infty)$. We give the generalization of this
to quasiperiodic operators in Theorem~\ref{thm:stabess}.

Kneser answered in \cite{kn}, the question when $0$ is an accumulation
point of eigenvalues below $0$. One has if
\be
\limsup_{x\to\infty} q(x) x^2 < - \frac{1}{4}
\ee
then $0$ is an accumulation point of eigenvalues, and if
\be
\liminf_{x\to\infty} q(x) x^2 > - \frac{1}{4}
\ee
then $0$ is not an accumulation point of eigenvalues. The periodic case
was answered by Rofe-Beketov (see here his recent monograph \cite{krb}).
The generalization to the quasiperiodic case is given in
Theorem~\ref{thm:osci}.

Once it is known that $0$ is an accumulation point of eigenvalues, it is natural
to ask how fast do the eigenvalues converge to $0$. This question was answered
by Kirsch-Simon in \cite{kisi}. To state their result let $N(\lam)$ be
the number of eigenvalues of $-\frac{d^2}{dx^2} + \frac{\mu}{x^2}$ below $\lam$, then
\be
N(\lam) = \frac{1}{4 \pi} \sqrt{\frac{\mu}{\mu_{crit}} - 1} \abs{\ln\abs{\lam}} (1 + o(1)),
\quad \lam\uparrow 0,\quad \mu_{crit} = - \frac{1}{4}.
\ee
For $-\frac{d^2}{dx^2} + \frac{\mu}{x^\gam}$, $0 < \gam < 2$, we have
\begin{align}
N(\lam) &= \frac{1}{\pi} \int_{\{x,\, q(x) < \lam\}} (\lam - q(x))^{1/2} dx (1 + o(1)),& \lam\uparrow 0 \\
\nn & = \frac{\sqrt{\mu/\mu_{crit}}}{\pi} \frac{1}{2 - \gam}
\left| \frac{\mu}{\lam} \right|^{(2-\gam)/2\gam} (1 + o(1)),& \lam\uparrow 0
\end{align}
see Theorem XIII.82 in \cite{rs4}.
\footnote{We obtain a factor $\frac{1}{2}$ different from \cite{kisi} in
the case $\gam =2$, since we are considering half line operators. This
factor does not arise for $0 < \gam < 2$, since the domain of
integration also shrinks.}
This result goes back to results in the
sixties, see the notes in \cite{rs4}.
The periodic case was answered by Schmidt \cite{schm} for $\gam = 2$.
We will answer this question in Theorem~\ref{thm:eigasym}.

Periodic operators have a spectrum made out of the union of
finitely or infinitely many bands. That is
\be
\sig_{ess}(-\frac{d^2}{dx^2} + q_0(x)) = [E_0, E_1] \cup [E_2, E_3] \cup \dots,
\quad E_j < E_{j+1},
\ee
for $q_0(x + p) = q_0(x)$, $p > 0$. Since, we now have several
boundary points of the spectrum, one can also ask
what happens at all, finitely many, \dots boundary points
of $\sig_{ess}(H_0)$. Rofe-Beketov gave the following answer to this question:
Only finitely many gaps can contain
infinitely many eigenvalues for critical perturbations ($q(x) = \mu/x^2$)
(see (6.145) in \cite{krb}).
We will treat this question in Theorem~\ref{thm:howmany}.

The organization of this paper is as follows. In Section~\ref{sec:quasi},
we will state the needed results about quasiperiodic Schr\"odinger operators.
In Section~\ref{sec:main}, we will state our main results. Most proofs are
easy enough to be directly stated. Only the eigenvalue asymptotics requires
more work and is stated in the following section. In Section~\ref{sec:outlinekam},
we give an outline of Eliasson's proof and derive some further estimates.
In Appendix~\ref{sec:relosc}, we will review relative oscillation theory,
followed by another short appendix on needed methods from the theory
of differential equations. 

\section{Quasiperiodic Operators}
\label{sec:quasi}

We will now recall the basic notations about quasiperiodic
Schr\"odinger operators.
Let $\mathbb{T}^d$ be the $d$-dimensional torus, where
$\mathbb{T} = \mathbb{R}/(2\pi\mathbb{Z})$. Let $Q:\mathbb{T}^d\to\R$ be a real analytic
function. We will consider the Schr\"odinger operator on $L^2(1,\infty)$ given by
\be
H_0 = - \frac{d^2}{dx^2} + q_0(x),\quad q_0(x) = Q(\omega x)
\ee
where $\omega\in\mathbb{T}^d$ is fixed.
We will assume that $\omega$ is a Diophantine number, that is there is some
$\tau > d - 1$, $\kappa>0$,
such that
\be\label{eq:conddiophane}
DC(\kappa, \tau): \quad \abs{\spr{\omega}{n}} \geq \frac{\kappa}{\abs{n}^\tau},\quad n\in\mathbb{Z}^d \backslash\{0\},
\ee
holds.

Recall the rotation number $\rho(E)$ from \cite{jm}.
Denote by $\vartheta(x,E)$ the Pr\"ufer angle of a solution
$u$ of $H_0 u = E u$. That is a continuous function of $x$
such that
\be
u(x) = r(x) \sin\vartheta(x,E),\quad u'(x) = r(x) \cos\vartheta(x,E),
\quad 0 \leq \vartheta(1,E) < \pi,
\ee
for some continuous function $r$. The rotation number $\rho(E)$ is now
introduced by
\be
\rho(E) = \lim_{x\to\infty} \frac{\vartheta(x,E)}{x}.
\ee
We remark that the integrated density of states $k(E)$
satisfies
\be
k(E) = \frac{1}{\pi} \rho(E).
\ee
Johnson and Moser showed 
\begin{theorem}\label{thm:jm}[\cite{jm}]
The spectrum $\sig(H_0)$ is given by
\be
\sig(H_0) = \{E,\,\rho(E)=\frac{1}{2} \spr{\omega}{n},\quad n\in\Z^d\}.
\ee
Furthermore $\rho$ is a continuous function, and constant outside the
spectrum.
\end{theorem}

Now we come to Eliasson's result.
Recall that we can rewrite the Schr\"odinger equation
$$
-u''(x) + Q( \omega x) u(x) = E u(x),
$$
as the first order system
\be\label{eq:matrixode}
 U'(x) = \begin{pmatrix} 0 & 1 \\ Q(\omega x) - E & 0 \end{pmatrix} U(x)
\ee
where $U(x) = \begin{pmatrix} u(x) \\ u'(x) \end{pmatrix}$.

\begin{theorem}\label{thm:formsol}[\cite{eli}]
There is an $E_0$, such that for $E = \frac{1}{2} \spr{\omega}{m} > E_0$
a boundary point of the spectrum of $H_0$, there is a function 
$Y: \mathbb{T}^d\to SL(2,\R)$ and $A\in sl(2,\R)$ with $A^2=0$ such
that
\be\label{eq:formsol}
X(x) = Y_1 Y(\frac{\omega}{2} x) e^{A x},\quad
Y_1=\frac{1}{2 \sqrt{E}}\begin{pmatrix}1&1\\-\sqrt{E}&\sqrt{E}\end{pmatrix}
\ee
is the fundamental solution of (\ref{eq:matrixode}).
Furthermore we have that for $\abs{m} \geq 2$
\begin{align}
\label{eq:estinorma2} \abs{A} &\leq c \abs{m}^{\frac{3}{2} \tau}\\
\label{eq:estinormY}
\abs{Y} &\leq C \log {\abs{m}},\quad \abs{\det(Y) - 1} \leq \frac{1}{2},
\end{align}
for constants $c$, $C$ independent of $m$, and
the spectrum of $H_0$ is purely absolutely continuous above $E_0$.
\end{theorem}

We will give an outline of Eliasson's proof in 
Section~\ref{sec:outlinekam}, and derive the additional
estimates there.
In fact Eliasson proved that (\ref{eq:formsol}) holds,
when $\rho(E)$ satisfies the next Diophantine condition
\be
\abs{\rho - \frac{\spr{n}{\omega}}{2}} \geq \frac{\ti{\kappa}}{\abs{n}^{\sigma}},\quad n\in\mathbb{Z}^d \backslash\{0\},
\quad \ti{\kappa} > 0,\,\sigma >0.
\ee
Eliasson also showed that the spectrum of $H_0$ will be
a Cantor set for generic functions $Q:\mathbb{T}^d \to \R$ in the $\abs{.}_s$
topology given by the norm
\be\label{eq:topanalytic}
| Q |_s = \sup_{\abs{\im(z)} < s}  \abs{Q(z)}.
\ee
Furthermore, we could replace $Q(\omega x)$ by $Q(\omega x + \theta)$
for any $\theta\in\mathbb{T}^d$ obtaining the same statement.

\section{Main Results}
\label{sec:main}

We are interested in decaying perturbations of the quasiperiodic 
operator $H_0$. That is for
some function $\Delta q$ consider the operator
\be
H_1 = -\frac{d^2}{dx^2} + q_1(x),\quad q_1(x) = q_0(x) + \Delta q(x),
\ee
for $q_0(x) = Q(\omega x)$ as described in Section~\ref{sec:quasi}.
We then have the next basic stability result of the essential spectrum.
\begin{theorem}\label{thm:stabess}
If $\Delta q \to 0$, then
\be
\sig_{ess}(H_1) = \sig_{ess}(H_0) = \R \backslash \bigcup_{n} G_n,
\ee
for open sets $G_n$. If $\Delta q$ is integrable, we have that
the spectrum of $H_1$ is purely absolutely continuous
above $E_0$.
\end{theorem}

\begin{proof}
The first part follows by Weyl's Theorem and Theorem~\ref{thm:jm}.
For the second part, note that by Theorem~\ref{thm:formsol}, $H_0$
has purely absolutely continuous spectrum above $E_0$, and by
Theorem 1.6. of \cite{kls} it is invariant under $L^1$ perturbations.
\end{proof}

It is conjectured in \cite{kls}, that there is still absolutely continuous
spectrum for $\Delta q \in L^2$, but it may not be pure.
This was shown for the free case in \cite{deikil} and for the periodic
one in \cite{kil}. See also the recent review in \cite{deki}.
If we write $G_n = (E_n^-, E_n^+)$ for the intervals of the last theorem,
and call them gaps. We call $E_n^-$ (resp. $E_n^+$) a lower (resp. upper) boundary
point of the spectrum. The next relative oscillation criterion follows.

\begin{theorem}\label{thm:osci}
Assume that $\Delta q \to 0$, and let $E$ be a boundary point above $E_0$
of the essential spectrum of $H_0$. Then there exists a constant
$K = K(E)$ such that $E$ is an accumulation point of eigenvalues of $H_1$
if
\be
\limsup_{x\to\infty} K \Delta q(x) x^2 < - \frac{1}{4}
\ee
and $E$ is not an accumulation point of eigenvalues if
\be
\liminf_{x\to\infty} K \Delta q(x) x^2 > - \frac{1}{4}.
\ee
Furthermore $K > 0$ (resp. $K < 0$) if $E$ is a upper (resp. lower) boundary
point.
\end{theorem}

\begin{proof}
Everything follows from Theorem~\ref{thm:main},
except for the existence of $K$. We have from (\ref{eq:formsol})
that $u_0(t) = U(\frac{\omega}{2} t)$ for a function $U:\mathbb{T}^d \to \R$.
We will show
\begin{align*}
\nn K &= \liminf_{l\to\infty} \limsup_{x\to\infty} \frac{1}{\ell} \int_x^{x+\ell} u_0(t)^2 dt  \\
&=  \limsup_{l\to\infty} \liminf_{x\to\infty} \frac{1}{\ell} \int_x^{x+\ell} u_0(t)^2 dt =
\int_{\mathbb{T}^d} U(z)^2 dz.
\end{align*}
Now note, that (\ref{eq:conddiophane}) implies that
the system $(\mathbb{T}^d, T_t, \mu)$, where $T_t = \frac{\omega}{2} t$ and
$\mu$ is the normalized Lebesgue measure is uniquely ergodic.
By Birkhoff's ergodic theorem, we have that
$$
\lim_{l\to\infty} \frac{1}{\ell} \int_x^{x+\ell} U(\frac{\omega}{2} t)^2 dt = \int_{\mathbb{T}^d} U(z)^2 dz.
$$
By unique ergodicity, we know that the limit is uniform in $x$. Hence, the result follows.
\end{proof}

We even have a whole scale of relative oscillation criteria.
To state this, we recall the iterated logarithm $\log_n(x)$ which is defined recursively via
$$
\log_0(x) = x, \qquad \log_n(x) = \log(\log_{n-1}(x)).
$$
Here we use the convention $\log(x)=\log|x|$ for negative values of $x$. Then
$\log_n(x)$ will be continuous for $x>\E_{n-1}$ and positive for $x>\E_n$, where
$\E_{-1}=-\infty$ and $\E_n=\E^{\E_{n-1}}$. Abbreviate further
$$
L_n(x) = \frac{1}{\log_{n+1}'(x)} = \prod_{j=0}^n \log_j(x),\qquad
\ti{Q}_n(x) = -\frac{1}{4 K} \sum_{j=0}^{n-1} \frac{1}{L_j(x)^2}.
$$
From Theorem 2.10. of \cite{ktep}.

\begin{theorem}
Assume the assumptions of the last theorem, and
that for some $n\in\N$
\be
\lim_{x\to\infty} L_{n-1}(x)^{-2} (\Delta q(x) - \ti{Q}_{n-1}(x)) = -\frac{1}{4 K}.
\ee
Then $E$ is an accumulation point of eigenvalues of $H_1$ if
\be
\limsup_{x\to\infty} K L_n(x)^2 (\Delta q(x)  - \ti{Q}_n(x)) < - \frac{1}{4}
\ee
and $E$ is not an accumulation point of eigenvalues if
\be
\liminf_{x\to\infty} K L_n(x)^2 (\Delta q(x) - \ti{Q}_n(x)) > - \frac{1}{4},
\ee
with the same $K$ as in the last theorem.
\end{theorem}

The next lemma gives us an estimate on $K$.

\begin{lemma}\label{lem:estiI}
The constant $K$ of Theorem~\ref{thm:osci}, satisfies
\be
\abs{K(E)} \leq \frac{C}{\abs{m}^{\ti{\tau}} \sqrt{E}},\quad 0<\ti{\tau} < \frac{3}{2} \tau
\ee
where $m \in \Z^d$ is such that $E \in \rho^{-1}(\frac{1}{2} \spr{\omega}{m})$.
\end{lemma}

\begin{proof}
From Theorem~\ref{thm:formsol}, we know the existence. We note that
$\det(Y_1) = 1$.
By (\ref{eq:estinorma2}), we have that $\abs{A} \leq c \abs{m}^{-\frac{3}{2} \tau}$, where
the $m$ is the one such that $\rho(E) = \frac{1}{2} \spr{\omega}{m}$. Hence,
we obtain that $\abs{\beta} \leq c \abs{m}^{-\frac{3}{2} \tau}$. Now
$$
\abs{K} \leq c \frac{\int_{\mathbb{T}} U(z) dz}{\abs{m}^{\frac{3}{2}\tau} \sqrt{E}}.
$$
The claim now follows by (\ref{eq:estinormY}).
\end{proof}

\begin{remark}
 One can hope that the estimate (\ref{eq:mucrit}) on $K(E)$ can be improved.
 It was shown in \cite{eli} that the matrix $A$ and then $\beta$ would
 satisfy the bound $\abs{\beta} \leq C \abs{E_+ - E_-}$ for some constant
 $C$. Then it was shown in \cite{mp}, that $\abs{E_+ - E_-} \leq c e^{-\gamma \abs{m}}$
 for some constants $c$ and $\gamma$. Hence one should expect
 $K(E)$ to decrease exponentially in $m$. Unfortunately, the estimate
 of \cite{mp} depends on further arithmetic properties of $m$. Hence,
 it is not clear if it holds at all band edges.
\end{remark}

For simplicity, we will now restrict our attention to perturbations of the form
\be
\Delta q (x) = \frac{\mu}{x^\gam},\quad \mu\neq  0, \quad \gam > 0.
\ee
We will denote the operator $H_0 +  \frac{\mu}{x^\gam}$ by $H^\gam_\mu$. Now, we come to
the question how many gaps above $E_0$ can contain infinitely many
eigenvalues. This question is a bit odder than the one for periodic operators,
since there are bounded intervals that contain infinitely many gaps.

Introduce $\mu_{crit}$ by
\be\label{eq:mucrit}
\mu_{crit}(E) = - \frac{1}{4 K(E)}.
\ee
Then $E > E_0$ is an accumulation point of eigenvalues of $H^2_\mu$ if
and only if $\mu/\mu_{crit} > 1$. For $H^\gam_\mu$, $0<\gam<2$,
this requirement is $\mu/\mu_{crit} > 0$. Now, we come to

\begin{theorem}\label{thm:howmany}
If $\gam > 2$, then no boundary point of $\sig(H_0)$ above $E_0$ is an accumulation
point of eigenvalues of $H^\gam_\mu = H_0 + \frac{\gam}{x^\alpha}$.
If $\gam < 2$, then if $\mu < 0$ (resp. $\mu > 0$),
then all upper (resp. lower) boundary points above $E_0$
are accumulation points of eigenvalues of $H^\gam_\mu$.

If $\gam = 2$, we can add infinitely many eigenvalues to each gap by choosing
$\mu$ large enough. However, for every value of $\mu$ only finitely many gaps
contain infinitely many eigenvalues of $H^\gam_\mu$.
\end{theorem}

\begin{proof}
The first claim follows from Theorem~\ref{thm:osci}. The second
claim follows from the last lemma and Theorem~\ref{thm:osci}.
\end{proof}

Now, we come to the eigenvalue asymptotics. Let $E$ be again a boundary
point of the spectrum of $H_0$.
Introduce, if the set $(\ti{E}, E) \cap \sig(H_0) = \emptyset$, by $N(\lam)$ the number
\be
N(\lam) = \tr(P_{(\ti{E}, \lam)} (H_\mu^\gam)),\quad \ti{E} < \lam < E
\ee
with the obvious modification for $(E,\ti{E}) \cap \sig(H_0) = \emptyset$.

\begin{theorem}\label{thm:eigasym}
Let $E$ be a boundary point of the spectrum of $H_0$, which is an accumulation
point of eigenvalues of $H^\gam_\mu$.  Then if $\gam = 2$
\be
N(\lam) = \frac{1}{4\pi} \sqrt{\frac{\mu}{\mu_{crit}} - 1} \cdot \abs{\log{\abs{E - \lam}}} \cdot (1 + o(1)),
\ee
and if $0 < \gam < 2$
\be
N(\lam) = \frac{1}{\pi} \frac{1}{2 - \gam} \sqrt{\frac{\mu}{\mu_{crit}}}
\left(\frac{\abs{\mu}}{\abs{E - \lam}}\right)^{(2 - \gam)/ 2\gam} \cdot (1 + o(1)).
\ee
where $N(\lam)$ is the number of eigenvalues near $E$.
\end{theorem}

We will give a proof in Section~\ref{sec:proofasym}. In difference to the
proof of \cite{schm}, our proof only uses the decay of the potential and
the behavior of the solution at the boundary point of the spectrum. In fact
everything carries over to general elliptic situations. That is,
where one has two solutions $u_0$, $u_1$ such that $u_0(x)$ and
$u_1(x) - x u_0(x)$ are bounded functions.

\begin{remark}
It was already shown in Corollary 6.6 in \cite{krb}
that $\mu_{crit}(E)$ has to diverge as $E \to \infty$. We also remark
that \cite{krb} develops a different approach to relative oscillation criteria
than was used in \cite{ktep}.
\end{remark}

\section{Proof of Theorem~\ref{thm:eigasym}}
\label{sec:proofasym}

We will now give explicit bounds on the spectral projections.

\begin{lemma}\label{lem:asymprue}
Let $\psi$ be a solution of
\be
\psi'(x) = - \Delta q(x) (u_0(x) \cos\psi(x) - v_0(x) \cos\psi(x))^2.
\ee
Then we have that
\be
\psi(x) = \left(\frac{1}{2} \sqrt{\frac{\mu}{\mu_{crit}} - 1} + o(1)\right) \log\abs{x}
\ee
if $\Delta q (x) = \mu/x^2$.
If $\Delta q(x) = \mu/x^\gam$, $0 < \gam < 2$,
\be
\psi(x) =  \left(\sqrt{\frac{\mu}{\mu_{crit}}} \frac{1}{2 - \gam} + o(1)\right) x^{1 - \gam/2}.
\ee
\end{lemma}

\begin{proof}
Use in (\ref{eq:diffvarphi}) $\alpha = x$, to obtain if
$\gam = 2$ the next equation
$$
\varphi' = \frac{1}{x} (\sin^2 \varphi + \cos \varphi \sin \varphi + 
\mu u_0^2 \cos^2 \varphi) + O(\frac{1}{x^2}),
$$
whose asymptotics can be evaluated with Lemma~\ref{lem:aver} and Lemma~\ref{lem:boundsol}.

In the case $0 < \gam < 2$, we choose $\alpha = 1/\sqrt{| \Delta q| K}$, then
also the $\sin\varphi\cos\varphi$ term becomes of lower order, hence we obtain by averaging
$$
\varphi'(x) = \sqrt{K \Delta q(x)}  + O(\Delta q + \frac{\Delta q'}{\Delta q})
$$
which implies the claim for $\Delta q$ of the particular form.
\end{proof}

Then, we have that

\begin{lemma}\label{lem:upperbound}
Let the Wronskian $W(u_1(E), u_0(E))$ have $n$ zeros on
$$
\{x,\,\forall y>x,\, \abs{\Delta q} \leq \abs{E - \lam}\}.
$$ 
Then we have that
\be
N(\lam) \leq n + 3.
\ee
\end{lemma}

\begin{proof}
Observe that by the comparison theorem for Wronskians, we have that
$W(u_1(\lam), u_0(\lam))$ can have at most one zero left of $x_n$.

Hence, we obtain
$$
\dim\Ran P_{(-\infty,\lam)}(H_1) \leq \#_{(1,x_n)}(u_1(\lam),u_0(\lam)) + 1
$$
Now, by the triangle inequality for Wronskians, we obtain
$\#_{(1,x_n)}(u_1(\lam),u_0(\lam)) \leq \#_{(1,x_n)}(u_1(\lam),u_0(E)) + 1$.
It, now suffices to note that $\#_{(1,x_n)}(u_1(\lam),u_0(E))$
is bounded by $\#_{(1,x_n)}(u_1(E),u_0(E)) + 1$ by using the comparison
theorem for Wronskians.
\end{proof}

Note, that the last two lemmas imply the next bound on the eigenvalues
if $\gam = 2$
\be\label{eq:upbounda2}
N(\lam) \leq \frac{1}{4\pi} \sqrt{\frac{\mu}{\mu_{crit}} - 1} \cdot \abs{\log{\abs{E - \lam}}} \cdot (1 + o(1)),
\ee
and if $0 < \gam < 2$
\be\label{eq:upboundal2}
N(\lam) \leq \frac{1}{\pi} \frac{1} {2 - \gam}\sqrt{\frac{\mu}{\mu_{crit}}} 
\left(\frac{\abs{\mu}}{\abs{E - \lam}}\right)^{(2 - \gam)/ 2\gam} \cdot (1 + o(1)).
\ee
The next lemma shows that we have equality.
Hence with it Theorem~\ref{thm:eigasym} is proven.

\begin{lemma}\label{lem:lowerbound2}
Let $0 < \delta < 1$ and $0 < \gam \leq 2$, then
if $\gam = 2$,
\be
N(\lam) \geq \frac{1}{4\pi} \sqrt{\frac{\mu}{\mu_{crit}} (1 - \delta) - 1} \cdot \abs{\log{\abs{E - \lam}}} \cdot (1 + o(1)),
\ee
and if $0 < \gam < 2$,
\be
N(\lam) \geq \frac{1}{\pi} \frac{1}{2 - \gam} \sqrt{\frac{\mu}{\mu_{crit}} (1 - \delta)} 
\left(\frac{\abs{\mu}}{\abs{E - \lam}}\right)^{(2 - \gam)/ 2\gam} \cdot (1 + o(1)).
\ee
\end{lemma}

\begin{proof}
Let $x_{max}$ be given by
$$
x_{max}(\lam) = \delta \left(\frac{\abs{\mu}}{\abs{E - \lam}}\right)^{1/\gam}.
$$
Let $\varphi_\lam(x)$ be a Pr\"ufer angle of $W(u_0(E), u_1(\lam))$.
By the triangle inequality for Wronskians, it is clear that
$\varphi_\lam$ is close to the Pr\"ufer angle of $W(u_0(\lam), u_1(\lam))$.
Now, for $x < x_{max}(\lam)$, we have that
$$
\varphi_\lam'(x) \geq \frac{\mu ( 1 - \delta)}{x^\gam}
(u_0 \cos\varphi_\lam(x) - v_0 \sin\varphi_\lam(x))^2.
$$
This is the same equation for all $\lam$. As $x\to\infty$,
the solution has the claimed asymptotics by using Lemma~\ref{lem:asymprue}.
Hence, the claim follows.
\end{proof}

\section{Outline of Eliasson's proof}
\label{sec:outlinekam}

We now give an outline of Eliasson's proof of reducibility
in \cite{eli}. The next lemma is an easy computation.

\begin{lemma}\label{lem:hypkam}
 The equation
 \be
 X'(x) = \begin{pmatrix} 0 & 1 \\ Q(\omega x) - E & 0 \end{pmatrix} X(x)
 \ee
 can be transformed by
 \be
 X_1(x) = Y_1^{-1} X(x),\quad Y_1=\begin{pmatrix}1&1\\-\sqrt{E}&\sqrt{E}\end{pmatrix}
 \ee
 to
 \be
 X_1'(x) = (A_1 + F_1(\omega x, \sqrt{E})) X_1(x),
 \ee
 where
 \be
 A_1 = \sqrt{E} J,\quad J = \begin{pmatrix} 0 & 1 \\ -1 & 0 \end{pmatrix},
 \quad F_1(z,\sqrt{E}) = \frac{Q(z)}{2 \sqrt{E}} \begin{pmatrix} -1 & -1 \\ 1 & 1 \end{pmatrix}.
 \ee
 Furthermore $A_1, F_1$ satisfy Hypothesis~\ref{hyp:indass}.
\end{lemma}

\begin{hypothesis}\label{hyp:indass}
Let $A_1 \in sl(2,\C)$ and $F_1:\mathbb{T}^d \to Mat(2,\C)$ satisfy
\begin{align}
\label{eq:assF1} \tr(\hat{F}_1(0)) & = 0, \\
\label{eq:assA1} \abs{A_1 - \sqrt{E} J} & < 2 & \\
\label{eq:assF2} \abs{F_1}_{r_1} & < \eps_1,
\end{align}
for some $\eps_1 > 0$, small.
\end{hypothesis}

We have now seen that we can reduce our system to one of the form
\be\label{eq:diffX1}
X_1'(x) = (A_1 + F_1(x))X_1(x),
\ee
where $F_1$ is small. This system although close to a constant
coefficient one cannot be solved explicitly.
However, we can reduce it to a system
\be\label{eq:diffX2}
X_2'(x) = (A_2 + F_2(x)) X_2(x)
\ee
where $F_2$ is smaller than $F_1$, as follows.
We will construct $A_2, F_2$, and a solution $Y_1$ of the
system
\be\label{eq:diffY}
Y_1'(x) = (A_1 + F_1) Y_1 - Y_1 (A_2 + F_2).
\ee
Then for $X_2$ a solution to (\ref{eq:diffX2}), we have that $Y_1 X_2$ will
solve (\ref{eq:diffX1}). Of course, we cannot hope that (\ref{eq:diffX2})
will be explicitly solvable, however we will be able to iterate the
above procedure to obtain better and better approximate solutions.

Since, we will require that $F_k \to 0$, and then
our final $X_\infty(x) = e^{x A}$. Here $A = \lim_{k\to\infty} A_k$.
So the final solution will be
\be
\prod_{k=1}^\infty Y_k(\frac{\omega}{2} x) e^{A x}
\ee
We will not attempt to solve (\ref{eq:diffY}) in this paper, and refer
to \cite{eli} for the details. However, we will draw further conclusions
from Eliasson's method to control our quantities.

Fix $0 < \eps_1 < 1$ sufficiently small. Fix $0 < \sigma < 1$, and let
\be
\eps_{j+1} = \eps_j^{1 + \sigma} = \eps_1 ^{(1 + \sigma)^{j}}.
\ee
Furthermore, assume that $r_1 > r_2 > r_3 > \dots$ is a decreasing sequence of positive numbers,
satisfying
\be
\frac{r_1}{2^{j+1}} \leq r_j - r_{j+1}.
\ee
$r_j$ will play the role of the neighborhood, where we suppose to have analyticity.
Introduce $N_j$ by
\be\label{eq:defNj}
N_j = \frac{2\sigma}{r_j - r_{j+1}} \log(\eps_j^{-1}) = \frac{2 \sigma (1+\sigma)^{j}}{r_j - r_{j+1}} \log(\eps_1^{-1})
\leq C (2 + 2 \sig)^{j},\quad C>0.
\ee
Furthermore, one also sees that $N_j \geq \ti{C} (1 + \sig)^{j}$ for some other constant $\ti{C}$.
Hence $N_j \to \infty$ as $j\to\infty$. Furthermore, we have
\be
\label{eq:estiepsbyN} \eps_j^\sig \leq \left(\frac{4\sig}{r_1(1 +\sig)}\log(\eps_1^{-1})(2+2\sig)^j\right)^{-4\tau}\leq N_j^{-4\tau}
\ee
if $\eps_1$ is small enough.

\begin{proposition}\label{prop:elilem12}
 Assume Hypothesis~\ref{hyp:indass} with $\eps_1$ small enough, then there are functions
 $Y_j : 2\mathbb{T}^d \to GL(2,\R)$,
 $A_j \in sl(2,\R)$, and $F_j: \mathbb{T}^d\to Mat(2,\R)$,
 for $j\geq 1$. Furthermore, there are numbers $m_j$ that
 satisfy
 \be\label{eq:defLambda}
 \eps_j^\sigma \leq \abs{2 \alpha_j - \spr{\omega}{m_j}} \leq 2 \eps_j^\sigma,\quad
 0 < \abs{m_j} \leq N_j,
 \ee 
 or $m_j = 0$ if (\ref{eq:defLambda}) cannot be satisfied. Here $\alpha_j$
 is the rotation number of $A_k$. Furthermore $A_j$, $F_j$, and $Y_j$ satisfy
 \begin{align}
 \label{eq:diffY1}
 \spr{Y_{j+1}'(x)}{\frac{\omega}{2}} &= (A_j + F_j(x)) Y_{j+1}(x) - Y_{j+1}(x) (A_{j+1}+F_{j+1}(x)),
 \end{align}
 \begin{align}
 \label{eq:estiY}
 \abs{\left(Y_{j+1}(.) - \exp\left(\frac{\spr{m_k}{.}}{\alpha_j} A_j\right)\right)} &
 \leq \eps_j^{1/2}, \\ 
 \label{eq:estiA1} \abs{\left(A_{j+1} - \left(1-\frac{\spr{\omega}{m_j}}{2\alpha_j}\right)A_j\right)} &
 \leq \eps_j^{2/3}, 
 \end{align}
 \begin{align}
 \label{eq:estiF} \tr(\hat{F}_{j+1}(0)) & = 0, \quad \abs{F_{j+1}}_{r_{j+1}} < \eps_{j+1},\\
 \label{eq:estiA3} \abs{A_{j+1}} &\leq 32 \abs{\alpha_{j+1}} N^\tau_{j+1},& \text{if } \abs{\alpha_{j+1}} &\geq \frac{1}{4} N_{j+1}^{-\tau}.
 \end{align}
 \end{proposition}
 
 \begin{proof}
 This is Lemma~1 and 2 in \cite{eli}.
 \end{proof}

\begin{remark}\label{rem:Elarge}
 The requirement of $\eps_1$ being small enough, will in fact determine
 our lower bound on allowed energies $E$. Since for $E > E_0$
 $$
 \abs{F_1}_{r_1} = \frac{C}{\sqrt{E}} < \frac{C}{\sqrt{E_0}} = \eps_1
 $$
 for some constant $C$. Hence by making $E_0$ large, we can make $\eps_1$ arbitrarily small.
\end{remark}

\begin{lemma}\label{lem:normAesti}
Assume that $Y_j$, $A_j$ and $F_j$ satisfy the conditions
given in Proposition~\ref{prop:elilem12}. 
If for all $j\leq k$, $m_j = 0$, then
\be
\abs{A_k - \lam J} < 3.
\ee
We furthermore obtain, if $K$ is the largest integer less than $k$ such that $m_K \neq 0$, that
\be\label{eq:estinorma}
\abs{A_k} \leq C \frac{1}{N_K^{3 \tau}} < 3,\quad k\geq K 
\ee
where $C$ doesn't depend on $K$.
\end{lemma}

\begin{proof}
By $m_k = 0$, we have that from (\ref{eq:defLambda})
$$
\abs{2\alpha_k - \spr{\omega}{n}}\geq \eps_k^\sig,\quad 0<\abs{n}\leq N_j.
$$
For $m_j =0$, $j=1,\dots,k$, we have by (\ref{eq:estiA1}) that
$$
\abs{A_k} \leq 2 + \eps_1^{2/3} + \dots + \eps_{k-1}^{3/2} < 3.
$$
This shows the first part.

For the second claim, let $l\leq k$ be maximal such that the $m_l\neq 0$.
Then, we obtain a bound on $\abs{A_j(\lam)}$ by
\begin{align*}
\abs{A_k} & \leq \eps_{l}^{2/3} + \dots + \eps_{k-1}^{3/2} + \abs{\left(1 - \frac{\spr{\omega}{m_l}}{2 \alpha_l} \right) A_l} \\
& \leq 2 \eps_l ^{2/3} + 2 \eps_l ^\sigma \abs{\frac{A_l}{2\alpha_l}} \leq 34 N_{l}^\tau \eps_l^\sigma,
\end{align*}
where we used (\ref{eq:defLambda}) in the middle and (\ref{eq:estiA3}) in the last step.
(\ref{eq:estinorma}) follows from (\ref{eq:estiepsbyN}). The last claim is evident.
\end{proof}

Let us now consider
$\ti{\rho} = \frac{1}{2}\sum_{k=1}^\infty \spr{m_k}{\omega} + \alpha$,
where $\alpha = \lim_{j\to\infty} \alpha_j$.
Furthermore,
$\rho_{j+1} = \frac{1}{2}\sum_{k=1}^j \spr{m_k}{\omega} + \alpha_{j+1}$,
Furthermore, we know that inside the gap $\alpha = 0$ from \cite{eli}.
We now obtain

\begin{lemma}
$\rho_{j+1} \to \ti{\rho}$ uniformly.
If $\rho$ is rational, $m_j = 0$ for $j$ large.
Furthemore,
\be\label{eq:summ}
\sum_{k, m_k \neq 0} m_k = m,
\ee
holds.
\end{lemma}

\begin{proof}
 The first two parts are Lemma~3 in \cite{eli}. The last part follows, since
 $\alpha\to 0$, and with $\ti{m} =  \sum_{k, m_k \neq 0} m_k$, one has
 $$
 0 = \ti{\rho} - \frac{1}{2} \spr{\omega}{m} = \frac{1}{2} \spr{\omega}{\ti{m} - m}.
 $$
 Hence $\ti{m} = m$ by the Diophantine condition.
\end{proof}

\begin{proof}[Proof of (\ref{eq:estinorma2})]
 We will now show how (\ref{eq:summ}) can be used to make the bound
 from (\ref{eq:estinorma}) only depend on $m$.
 By definition $\abs{m_k} \leq N_k$, we have by (\ref{eq:defNj})
 $$
 \abs{m} \leq \sum_{k, m_k\neq 0} \abs{m_k} \leq \sum_{k=1}^K N_k\leq C (2\sig + 2)^{K+1}.
 $$
 Hence $K \geq \frac{\log\abs{m}}{\log(2+2\sig)} - C$ and by (\ref{eq:defNj}) 
 $N_K \geq C \sqrt{\abs{m}}$
 and then (\ref{eq:estinorma}) implies the claim, since it holds for all large $k$.
\end{proof}

We have that
\begin{lemma}
If $m_j = 0$ for $j$ large,
we have that $\prod Y_j$ converges to some $Y$ uniformly on compact subsets.
Furthermore $A_j \to A$ and $F_j \to 0$. Furthermore (\ref{eq:estinormY}) holds.
\end{lemma}

\begin{proof}
 Since $r_j$ is decreasing and positive, it has a limit $r_0\geq 0$.
 By (\ref{eq:estiF}), we have that $\abs{F_j}_{r_j}\to 0$. Since $m_j = 0$
 for large $j$, we have that $\abs{A_{j+1} - A_j} \leq \eps_j^{2/3}$ from (\ref{eq:estiA1}).
 Hence, $A_j\to A$, since $\sum_{j=N}^\infty \eps_j^{2/3} < \infty$.
 
 By (\ref{eq:estiY}), we have that $\abs{Y_j - \mathbb{I}}\leq \eps_j^{1/2}$, if $m_j = 0$,
 which implies $\prod Y_j \to Y$ by a similar argument. 
 If $m_j\neq 0$, we have that
 $$
 \abs{Y_{j+1}(\frac{\omega}{2} t) - \mathbb{I}}\leq 
 \eps_j^{1/2} + \abs{\cos(\frac{\spr{m_j}{t}}{2})\mathbb{I} + \sin(\frac{\spr{m_j}{t}}{2})\frac{A_j}{\alpha_j} - \mathbb{I}},
 $$
 where the last term $\leq 3$. Since, we can bound the number of these terms by
 $\log{m}$, we obtain the claim.
 By (\ref{eq:estiY}), we have that $Y_{j+1} - \mathbb{I}$, if $m_j = 0$, resp.
 $\exp(- \spr{m_j}{t} A_j/\alpha_j) Y_{j+1} - \mathbb{I}$ are bounded by $\eps_{j}^{1/2}$.
 Hence, we can bound $\abs{\det(Y_{j+1}) - 1} \leq \eps_j$, from which the estimate
 on the determinant follows.
\end{proof}

\appendix

\section{Relative Oscillation Theory}\label{sec:relosc}

As introduced in \cite{kt1}, relative oscillation theory is a tool
to compute the difference of spectra of two different Schr\"odinger
operators. Let $q_0, q_1 \in L^1_{loc}$ and
\be
H_j = -\frac{d^2}{dx^2} + q_j,\quad j=0,1
\ee
be self-adjoint Schr\"odinger operators on $L^2(1,\infty)$.
Introduce $\Delta q = q_1 - q_0$, which we will assume to be
sign-definite. Denote by $\#(u_0, u_1)$ the number of
zeros of the Wronskian $W(u_0,u_1) = u_0 u_1' - u_0' u_1$ on $(1,\infty)$,
for solutions $\tau_j u_j = \lam_j u_j$. Let $\psi_{j,-}(\lam)$ be the solution
of $\tau_j \psi_{j,-}(\lam) = \lam \psi_{j,-}(\lam)$, which obeys the
boundary condition at $1$ (e.g. $\psi_{j,-}(\lam)(1) = 0$). Similarly
let $\psi_{j,+}(\lam)$ be the solution satisfying $\psi_{j,+}(\lam) \in L^2(1,\infty)$.
Then \cite{kt1} tells us:

\begin{theorem}\label{thm:mainosc}
 Assume that $[\lam_0,\lam_1] \cap \sig_{ess}(H_0) = \emptyset$. Then,
 we have that
 \begin{align}
 \nn \tr  P_{[\lam_0,\lam_1)}(H_1) &- \tr P_{(\lam_0,\lam_1]}(H_0) \\
 &= \begin{cases}
 (\#(\psi_{1,\pm}(\lam_1),\psi_{0,\mp}(\lam_1)-\#(\psi_{1,\pm}(\lam_0),\psi_{0,\mp}(\lam_0)),
 &  \Delta q < 0 \\         
 - (\#(\psi_{1,\pm}(\lam_1),\psi_{0,\mp}(\lam_1)-\#(\psi_{1,\pm}(\lam_0),\psi_{0,\mp}(\lam_0)),
 & \Delta q > 0
 \end{cases}
 \end{align}
 Here $\tr P_{[\lam_0,\lam_1)}(H_1)$ denotes the number of eigenvalues
 of $H_1$ in $[\lam_0, \lam_1)$.
\end{theorem}

Since one has the next triangle inequality for Wronskians
\be \label{eq:wtriang}
{\#}(u_0,u_1) + {\#}(u_1,u_2) - 1 \leq {\#}(u_0,u_2) \leq
{\#}(u_0,u_1) + {\#}(u_1,u_2) + 1,
\ee
one can replace $\psi_{j,\pm}(\lam)$ by any other solution of
$\tau_j u = \lam u$, up to a finite error. 
We furthermore remark that the next two comparison theorems
hold. The first one is found in \cite{kt2}.

\begin{theorem}[Sturm's Comparison theorem]\label{thm:sturclaform}
Let $q_0-q_1 > 0$, and $H_j u_j = 0$, $j=0,1$.
Then between any two zeros of $u_0$ or $W(u_0, u_1)$,
there is a zero of $u_1$.

Similarly, between two zeros of $u_1$, which are not at the same time zeros of $u_0$,
there is at least one zero of $u_0$ or $W(u_0,u_1)$.
\end{theorem}

The next result is found in \cite{kt1}.

\begin{theorem}[Comparison theorem for Wronskians]\label{thmscw}
Suppose $u_j$ satisfies $\tau_j u_j = \lam_j u_j$, $j=0,1,2$,
where $\lam_0 r -q_0 \le \lam_1 r - q_1 \le \lam_2 r - q_2$.

If $c<d$ are two zeros of $W_x(u_0,u_1)$ such that $W_x(u_0,u_1)$ does not
vanish identically, then there is at least one sign flip of $W_x(u_0,u_2)$ in $(c,d)$.
Similarly, if $c<d$ are two zeros of $W_x(u_1,u_2)$ such that $W_x(u_1,u_2)$ does not
vanish identically, then there is at least one sign flip of $W_x(u_0,u_2)$ in $(c,d)$.
\end{theorem}

We call $H_1$ relatively oscillatory with respect to $H_0$ at $E$ if for any solutions
of $H_j u_j(E) = E u_j(E)$, $j=0,1$, we have that $\#(u_0(E),u_1(E))$ is infinite.
Otherwise we call it relatively nonoscillatory.
Now, we come to relative oscillation criteria. 

\begin{lemma}\label{lem:nonoscingap}
Let $\lim_{x\to \infty} \Delta q (x) = 0$.
Then $\sig_{ess}(H_0)=\sig_{ess}(H_1)$ and
$H_1$ is relatively nonoscillatory with respect
to $H_0$ at $E \in \R \backslash\sigma_{ess}(H_0)$.
\end{lemma}

By Theorem~\ref{thm:mainosc}, this is equivalent if $E$ is a boundary
point of the essential spectrum of $H_0$, to $E$ being an accumulation point
of eigenvalues of $H_1$. In order to state a relative oscillation
criterion at a boundary point of the spectrum, some preparations are needed.

\begin{definition}\label{def:regedge}
A boundary point $E$ of the essential spectrum of $H_0$ will be called admissible if there
is a minimal solution $u_0$ of $(\tau_0-E) u_0= 0$ and a second linearly independent
solution $v_0$ with $W(u_0,v_0)=1$ such that
$$
\begin{pmatrix} u_0 \\ p_0 u_0' \end{pmatrix} = O(\alpha), \quad
\begin{pmatrix} v_0 \\ p_0 v_0' \end{pmatrix} - \beta\begin{pmatrix} u_0 \\ p_0 u_0' \end{pmatrix} =
o(\alpha \beta)
$$
for some weight functions $\alpha>0$, $\beta\lessgtr 0$, where $\beta$ is absolutely continuous such
that $\rho=\frac{\beta'}{\beta}>0$ satisfies $\rho(x)=o(1)$ and $\frac{1}{\ell} \int_0^\ell \left|\rho(x+t) -\rho(x) \right| dt = o(\rho(x))$.
\end{definition}

It is shown in Lemmas~4.2 and 4.3 of \cite{ktep},
that there exists a Pr\"ufer angle $\psi$ for $W(u_0, u_1)$ such that it obeys
\be
\psi'(x) = - \Delta q(x) (u_0(x) \cos(\psi(x)) - v_0(x) \sin(\psi(x)))^2.
\ee
Through the transformation $\cot\psi = \alpha \cot\varphi + \beta$,
this can then be transformed to (see Lemma~4.6 of \cite{ktep})
\begin{align}\label{eq:diffvarphi}
\varphi' &= \frac{\alpha'}{\alpha} \sin\varphi \cos\varphi + \frac{\beta'}{\alpha} \sin^2 \varphi
-\Delta q \cdot \alpha u_0^2 \cos^2 \varphi \\
\nn &+ O(\Delta q) + O(\Delta q/\alpha).
\end{align}

Through an application of the methods of Appendix~\ref{sec:ode}, 
one comes to the main result of \cite{ktep}.

\begin{theorem}\label{thm:main}
Suppose $E$ is an admissible boundary point of the essential spectrum of $\tau_0$,
with $u_0$, $v_0$ and $\alpha$, $\beta$ as in Definition~\ref{def:regedge}.
Furthermore, suppose that we have
\be
\Delta q = O\big(\frac{\beta'}{\alpha^2 \beta^2}\big).
\ee
Then $\tau_1-E$ is relatively oscillatory with respect to $\tau_0-E$ if
\be
\inf_{\ell>0} \limsup_{x\to b} \frac{1}{\ell} \int_x^{x+\ell} \frac{\beta(t)^2}{\beta'(t)}
u_0(t)^2 \Delta q(t) dt < -\frac{1}{4}
\ee
and relatively nonoscillatory with respect to $\tau_0-E$ if
\be
\sup_{\ell>0} \liminf_{x\to b} \frac{1}{\ell} \int_x^{x+\ell} \frac{\beta(t)^2}{\beta'(t)}
u_0(t)^2 \Delta q(t) dt > -\frac{1}{4}.
\ee
\end{theorem}

We finish this section with a closing remark.

\begin{remark}
 The requirement made that $\Delta q$ is of definite sign is not necessary.
 However, a general theory requires a more difficult definition of $\#(u_0,u_1)$.
 We refer the interested reader to \cite{kt1} for details.
\end{remark}

\section{Averaging ordinary differential equations}
\label{sec:ode}

In this section
we collect the required results for these ordinary differential equations.
Proofs and further references can be found in \cite{ktep}.

\begin{lemma}\label{lem:boundsol}
Suppose $\rho(x)>0$ (or $\rho(x)<0$) is not integrable near $b$. Then the equation
\be
\varphi'(x) = \rho(x) \bigg(A \sin^2\varphi(x) + \cos\varphi(x)\sin\varphi(x)
+ B \cos^2\varphi(x)\bigg) + o(\rho(x))
\ee
has only unbounded solution if $4 A B > 1$ and only bounded solution if $4 A B < 1$.
In the unbounded case we have
\be
\varphi(x) = \left(\frac{\sgn(A)}{2} \sqrt{4 A B-1} + o(1) \right) \int^x \rho(t) dt.
\ee
\end{lemma}

\noindent
In addition, we also need to look at averages:
Let $\ell > 0$, and denote by
\be
\ol{g}(x) = \frac{1}{\ell} \int_x^{x+\ell} g(t) dt.
\ee
the average of $g$ over an interval of length $\ell$.

\begin{lemma}\label{lem:aver}
Let $\varphi$ obey the equation
\be
\varphi'(x) = \rho(x) f(x) + o(\rho(x)),
\ee
where $f(x)$ is bounded. If
\be\label{condrho}
\frac{1}{\ell} \int_0^\ell \left|\rho(x+t) -\rho(x) \right| dt = o(\rho(x))
\ee
then
\be
\ol{\varphi}'(x) = \rho(x) \ol{f}(x) + o(\rho(x))
\ee

Moreover, suppose $\rho(x)=o(1)$. If $f(x)= A(x) g(\varphi(x))$, where $A(x)$ is bounded
and $g(x)$ is bounded and Lipschitz continuous, then
\be
\ol{f}(x)=  \ol{A}(x) g(\ol{\varphi}) + o(1).
\ee
\end{lemma}

\noindent
Condition (\ref{condrho}) is a strong version of saying that $\ol{\rho}(x)=\rho(x)(1+o(1))$ (it is
equivalent to the latter if $\rho$ is monotone). It will be typically fulfilled if $\rho$
decreases (or increases) polynomially (but not exponentially). For example, the condition holds if
$\sup_{t\in[0,1]} \frac{\rho'(x+t)}{\rho(x)} \to 0$.

Furthermore, note that if $A(x)$ has a limit, $A(x)=A_0+o(1)$, then $\ol{A}(x)$ can be replaced by
the limit $A_0$ and we have the next result

\begin{corollary}\label{cor:aver}
Let $\varphi$ obey the equation
\be
\varphi' = \rho \bigg(A \sin^2(\varphi) + \sin(\varphi) \cos(\varphi)
+ B \cos^2(\varphi)\bigg) + o(\rho)
\ee
with $A, B$ bounded functions and assume that $\rho=o(1)$ satisfies (\ref{condrho}).
Then the averaged function $\ol{\varphi}$ obeys the equation
\be
\ol{\varphi}' = \rho \bigg(\ol{A} \sin^2(\ol{\varphi}) + \sin(\ol{\varphi}) \cos(\ol{\varphi})
+ \ol{B} \cos^2(\ol{\varphi})\bigg) + o(\rho).
\ee
\end{corollary}

\noindent
Note that in this case $\varphi$ is bounded (above/below) if and only if $\ol{\varphi}$ is bounded
(above/below).

\subsection*{Acknowledgments}

I thank D.Damanik for bringing \cite{eli}, \cite{deki} to my attention,
and many helpful discussions.

\end{document}